\documentclass{svproc}
\usepackage{url}
\usepackage{amsmath}
\usepackage{amsfonts,amssymb}
\usepackage{graphicx}
\usepackage{tikz}
\usepackage{tikz-3dplot}
\usepackage{bigints}
\usetikzlibrary{patterns}
\usepackage{caption}
\usepackage{subcaption}
\usepackage{cite}

\allowdisplaybreaks

\begin{document}

\mainmatter              
\title{Discrete Generating Series and Linear Difference Equations}
\titlerunning{Discrete Generating Series}  
%
\author{V.~S.~Alekseev\inst{1} \and
T.~Cuchta\inst{2} \and A.~P.~Lyapin\inst{1}}
\authorrunning{Vitaliy Alekseev et al.} 
%
\tocauthor{Vitaliy Alekseev, Tom Cuchta, and Alexander Lyapin}
\institute{Siberian Federal University, Krasnoyarsk 660041, Russia
\and 
Marshall University, Huntington, WV USA 25755}

\maketitle              

\begin{abstract}
We define discrete generating series for arbitrary functions \( f \colon \mathbb{Z}^n \rightarrow \mathbb{C} \) and derive functional relations that these series satisfy. For linear difference equations with constant coefficients, we establish explicit functional equations linking the generating series to the initial data, and for equations with polynomial coefficients, we introduce an analogue of Stanley's \( D \)-finiteness criterion, proving that a discrete generating series is \( D \)-finite if and only if the corresponding sequence is polynomially recursive. The framework is further generalized to multidimensional settings, where we investigate the interplay between discrete generating series and solutions to Cauchy problems for difference equations. Key structural properties are uncovered through the introduction of polynomial shift operators and projection techniques. The theory is illustrated with concrete examples, including the tribonacci recurrence and Schröder's second problem. 
\end{abstract}

\section{Introduction}

In this paper, we develop a theory of discrete generating series associated with arbitrary functions \( f \colon \mathbb{Z}^n \to \mathbb{C} \), focusing on their connection to linear difference equations. We begin by defining such series for $n=1$ and derive the functional relations they satisfy, first for one-dimensional difference equations with constant coefficients and in turn for polynomial coefficients. Afterwards, we generalize the theory to partial difference equations in both cases. Our approach builds upon classical work on generating functions of linear recurrences, including the Stanley hierarchy (rational, algebraic, $D$-finite) \cite{stanley1980, Stanley1990, LeinartasYakovleva2018b}. The classification of generating functions and recursive sequences has been deepened in recent studies \cite{BousquetMelouPetkovsek2000, LyapinCuchta2022, Merlini2024}, particularly for $C$-finite and polynomially recursive ($p$-recursive) sequences.

Difference equations with polynomial coefficients have proven effective in enumerative combinatorics, especially in the analysis of restricted lattice paths \cite{Chandragiri2019, Chandragiri2023} and discrete special functions theory \cite{assaf2023discrete,BohnerCuchta2017}, structural and geometric aspects such as the amoeba of characteristic polynomials \cite{Krasikov2023}, and multidimensional versions of Poincaré’s theorem \cite{LeinartasPassareTsikh2008} offer additional insights into solution spaces \cite{Leinartas2004a}. Efficient techniques for deriving explicit formulas for generating series coefficients, particularly
in the context of the Aztec diamond and permutations with cycles, have been
developed in \cite{Kruchinin2021, Kruchinin2023}. Further developments on infinite linear difference operators and symbolic representations are presented in \cite{Abramov2011, Abramov2021, Abramov2015} and solutions of summation problems for polynomial difference operators and descriptions of the associated polynomial solution spaces were addressed in \cite{Grigoriev2022, LeinartasShishkina2022}.  We note that almost periodic and almost automorphic solutions to multidimensional difference equations considered in \cite{Kostic2023}.

In the context of discrete generating series, we introduce a notion of D-finiteness and prove a one-dimensional analogue of Stanley’s theorem about differential operators. We then generalize the theory to the multidimensional case for functions \( f : \mathbb{Z}^n \to \mathbb{C} \), where we introduce multidimensional generating series, polynomial shift operators, and projection techniques. These tools enable the formulation and solution of the Cauchy problems for linear difference equations with constant or polynomial coefficients.

The structure of the paper is as follows. In Section~\ref{sec:onedim}, we study the one-dimensional case, derive functional equations for discrete generating series, and prove a D-finiteness criterion for solutions of linear difference equations. Section~\ref{sec:multidim} extends the theory to multidimensional series, introducing the necessary operators and techniques to address the Cauchy problem. In Section~\ref{sec:examples}, we illustrate the theoretical framework with several examples, including the tribonacci recurrence and Schröder’s second problem. Finally, Section~\ref{sec:conclusion} summarizes our contributions and outlines future directions, including potential extensions to time scales and applications in combinatorics and signal processing.

\section{Discrete Generating Series: One-Dimensional Case}\label{sec:onedim}
\smallskip

Let $\mathbb N, \mathbb Z, \mathbb Z_\geqslant, \mathbb C$ denote the sets of natural numbers, integers, non-negative integers, and complex numbers, respectively. Let $f\colon \mathbb Z \to \mathbb C$, $z \in \mathbb C$, $n, \ell \in \mathbb Z_\geqslant$, and let $z^{\underline{n}} = z(z-1)\cdots (z-n+1)$ denote the discrete monomials (also known as ``falling factorial of $z$"). Let $(a)_k = a (a+1) \cdots (a+k-1)$ denote the Pochhammer symbol.

\subsection{Discrete Generating Series for Solutions of Linear Difference Equations with Constant Coefficients}

The simple connection between generating series and discrete generating series was given in \cite{Alekseev2023a}.

Define the forward shift operator $\delta$ by $\delta f(z)=f(z+1)$ and let the backwards shift operator be $\rho$, which obeys $\rho f(z)=f(z-1)$. Let $c_k \in \mathbb{C}$ be constants for $k=0, \ldots, r$. We define the polynomial difference operator as $P(\delta) = \displaystyle\sum_{k=0}^{r} c_k \delta^k$ and for $0 \leqslant \tau \leqslant r$, its truncation $\mathcal P_\tau(\delta) = \displaystyle\sum_{k=0}^{\tau} c_k \delta^k$.
\begin{definition}
A homogeneous linear difference equation with constant coefficients is a relation of the form
\begin{align}\label{difeqconst}
    P(\delta)f(x) := \sum_{k=0}^{r} c_k f(x+k) = 0.
\end{align}
\end{definition}
Equation \eqref{difeqconst} is often equipped with the Cauchy (or initial) data 
\begin{align}\label{initdata}
f(x) = \varphi(x), \quad x = 0, \ldots, r-1.
\end{align}
Solving the Cauchy problem \eqref{difeqconst}--\eqref{initdata} means finding a function $f\colon \mathbb Z \to \mathbb C$ that satisfies \eqref{difeqconst} and \eqref{initdata}. Such equations arise in a wide range of problems in enumerative combinatorics, such as lattice path problems or ballot problems (see an example in \cite{Nekrasova2014}), and in the theory of digital recursive filters \cite{Dudgeon1983}. The correctness of Cauchy problem for a polynomial difference operator is studied in  \cite{ApanovichLeinartas2017}.

\begin{definition}
The discrete generating series of $f$ is defined as the formal sum
\begin{equation}
F(\xi,\ell;z) = \sum_{x=0}^\infty f(x) \xi^x z^{\underline{\ell x}},\label{def:disgenseries}
\end{equation}
and for $s \in \mathbb Z_\geqslant$, consider we also consider truncations of \eqref{def:disgenseries}
\begin{equation*}
    F_s(\xi,\ell;z) = \sum_{x=0}^{s} f(x) \xi^x z^{\underline{\ell x}}.
\end{equation*}
\end{definition}
To study \eqref{def:disgenseries}, we need polynomial discrete shift operator
\begin{equation*}
    \mathcal R(\xi, \ell; z) = \sum_{k=0}^{r} c_{k}\xi^k z^{\underline{\ell k}}\rho^{\ell k}.
\end{equation*}
We now show that applying $\mathcal{R}$ to $F$ gives a polynomial in the discrete monomials.
\begin{theorem}\label{thm1}
The discrete generating series $F$ associated to \eqref{difeqconst}--\eqref{initdata} satisfies the functional relations
\begin{align*}
   \mathcal R(\xi, \ell; z) F(\xi,\ell;z) 
   &= \sum_{k=0}^{r-1} c_k  \xi^k z^{\underline{\ell k}} \rho^{\ell k} F_{r-k-1}(\xi,\ell;z) \\
   &= \sum_{x = 0}^{m-1} P(\delta^{-1}) \varphi(x) \xi^x z^{\underline{\ell x}} \\
   &= \sum_{x = 0}^{m-1} \mathcal P_{m-x}(\xi,\ell;z) \varphi(x) z^{\underline{\ell x}}.
\end{align*}
\end{theorem}
\begin{proof}
Multiplying both sides of \eqref{difeqconst} by $\xi^{x+r} z^{\underline{\ell (x+r)}}$ and summing over all $x \geqslant 0$
completes the proof after routine manipulations, see \cite{Alekseev2023a}. 
\end{proof}
Different groupings for the more general multidimensional case of generating series for Theorem~\ref{thm1} were presented in \cite{Leinartas2009}.

\subsection{$D$-Finite Discrete Generating Series}
In this section, we consider discrete generating series for solutions of linear difference equations with polynomial coefficients and prove an analogue of R. Stanley's theorem that derives conditions for the $D$-finiteness of discrete generating series \cite{stanley1980}.

Let $\Delta f(z) = f(z+1) - f(z)$ denote the forward difference operator. Then $\Delta z^{\underline{x}} = x z^{\underline{x-1}}$. Thus, the operator $\Delta$ and the discrete monomial are a discrete analogue of the power rule for classical derivatives. Further,
$$
\Delta F(\xi;\ell;z) = \Delta \sum\limits_{x=0}^\infty f(x) \xi^x z^{\underline{\ell x}} = \sum\limits_{x=0}^\infty \ell x f(x) \xi^x z^{\underline{\ell x - 1}}.
$$
Define the operator $\theta = \ell^{-1} z \rho \Delta$, see \cite{BohnerCuchta2018, CuchtaGrowWintz2023}. Consequently, 
\begin{equation} 
\theta z^{\underline{\ell x}}=xz^{\underline{\ell x}}.\label{thetaofdiscretemonomial}
\end{equation}
Let $p(x) = a_0 + a_1 x+ a_2 x^2 + ... + a_n x^n$ be a polynomial of degree $n$, and $p(\theta) = a_0 + a_1 \theta + a_2 \theta^2 + ... + a_n \theta^n$ be a polynomial difference operator. With these definitions, we know from \cite{Alekseev2023a} that
\[\theta^n F(\xi;\ell;z) = \sum\limits_{x=0}^\infty x^n f(x) \xi^x z^{\underline{\ell x}},\]
and 
\begin{equation} 
\sum\limits_{x=0}^\infty p(x) f(x) \xi^x z^{\underline{\ell x}} = p(\theta) \sum\limits_{x=0}^\infty f(x) \xi^x z^{\underline{\ell x}}.\label{pullpolynomialoutofsum}
\end{equation}
Now we prove a simple lemma.
\begin{lemma}\label{lemma3}
The following relation holds for all $x=1,2,3,\ldots$: 
\[(\theta - x) z^{\underline{\ell x}} = 0.\]
\end{lemma}

\begin{proof}
Using \eqref{thetaofdiscretemonomial}, compute
\[(\theta - x) z^{\underline{\ell x}} = \theta z^{\underline{\ell x}} - x z^{\underline{\ell x}} = x z^{\underline{\ell x}} - x z^{\underline{\ell x}} = 0,\]
which completes the proof.
\end{proof}
We now introduce a generalization of the difference equation \eqref{difeqconst}.
\begin{definition}
Let $\{p_k(x)\}_{k=0}^{r}$ be a set of polynomials, where $p_r(x)$ is not identically zero. A homogeneous linear difference equation with polynomial coefficients is a relation of the form
\begin{align}\label{difeqpol}
    \sum_{k=0}^{r} p_k(x) f(x+k) = 0.
\end{align}
Solving the Cauchy problem \eqref{initdata}, \eqref{difeqpol}  means finding a function $f\colon \mathbb Z \to \mathbb C$ that satisfies \eqref{initdata} and \eqref{difeqpol}.
\end{definition}

In \cite{stanley1980}, sequences $f$ for which there exist a polynomial sequence $p_k$ that satisfy \eqref{difeqpol} are referred to as polynomially recursive ($p$-recursive). The same work introduces the notion of $D$-finite power series, where the generating series $F(z)$ is termed $D$-finite if there exists a collection of polynomials $P_0(z)$, $P_1(z), \ldots$, $P_r(z)$, with at least one being non-zero, such that
\begin{align}\label{analog}
\left( P_r(z) \frac{\partial\,^r}{\partial z^r} + \cdots + P_1(z) \frac{\partial}{\partial z} + P_0(z)\right) F(z) = 0.
\end{align}
It was proven that the generating series $F(z)$ is $D$-finite if and only if the sequences $f(x)$ is $p$-recursive \cite[Theorem 1.5]{stanley1980}. The analogous notion of $D$-finiteness for discrete generating series resembles \eqref{analog}.
\begin{definition}
The discrete generating series $F(\xi,\ell;z)$ is called $D$-finite if there exists a collection of polynomials $$P_j(\xi,\ell;z) = \sum_{i=0}^{k_j} d_{i,j} \xi^{i} z^{\underline{i \ell}} \rho^{i \ell}, d_{i, j} \in \mathbb C, j = 0, \ldots, r,$$ with at least one being non-zero such that $F(\xi,\ell;z)$ solves the equation
\begin{align}\label{discrete}
\Big( P_r(\xi,\ell;z) \theta^r + \cdots + P_1(\xi,\ell;z) \theta + P_0(\xi,\ell;z)\Big) F(\xi,\ell;z) = 0.
\end{align}
\end{definition}
Now we prove an analogue of \cite[Theorem 1.5]{stanley1980}.
\begin{theorem}
The discrete generating series $F(\xi, \ell; z)$ is $D$-finite if and only if the sequence $f(x)$ is $p$-recursive.
\end{theorem}
\begin{proof}
First assume that the series $F(\xi,\ell;z)$ \eqref{def:disgenseries} associated to \eqref{difeqpol} is $D$-finite. Multiply both sides of \eqref{difeqpol} by $\xi^{x+r} z^{\underline{\ell (x+r)}}$ and sum over $x \geqslant 0$ to obtain
 \begin{align*}
   \sum\limits_{x=0}^\infty \sum\limits_{k=0}^{r} p_k(x)f(x+k) \xi^{x+r} z^{\underline{\ell (x+r)}} = 0.
\end{align*}
Commuting the summations, transforming the left-hand side, and routine manipulations yields
\[\begin{array}{ll}
   0&=\displaystyle\sum\limits_{x=0}^\infty \sum\limits_{k=0}^{r} p_k(x)f(x+k) \xi^{x+r} z^{\underline{\ell (x+r)}} \\
   &=\displaystyle\sum\limits_{k=0}^{r} \sum\limits_{x=k}^\infty p_k(x-k)f(x) \xi^{x+r-k} z^{\underline{\ell (x+r-k)}} \\
   &=\displaystyle\sum\limits_{k=0}^{r} \xi^{r-k} z^{\underline{\ell(r-k)}} \rho^{\ell(r-k)} \sum\limits_{x=k}^\infty p_k(x-k)f(x) \xi^{x} z^{\underline{\ell x}}. \\ 
\end{array}\]
Now let $q_k$ be the polynomial $q_k(x)=p_k(x-k)$. Thus we have derived
\[0=\displaystyle\sum\limits_{k=0}^{r} \xi^{r-k} z^{\underline{\ell(r-k)}} \rho^{\ell(r-k)} \sum\limits_{x=k}^\infty q_k(x)f(x) \xi^{x} z^{\underline{\ell x}},\]
and we use \eqref{def:disgenseries} and \eqref{pullpolynomialoutofsum} to get
\[\begin{array}{ll}
0&=\displaystyle\sum\limits_{k=0}^{r} \xi^{r-k} z^{\underline{\ell(r-k)}} \rho^{\ell(r-k)} q_k(\theta)\sum\limits_{x=k}^\infty f(x) \xi^{x} z^{\underline{\ell x}} \\ 
&=\displaystyle\sum\limits_{k=0}^{r} \xi^{r-k} z^{\underline{\ell(r-k)}} \rho^{\ell(r-k)} \Big[ q_k(\theta) F(\xi,\ell;z) -\displaystyle\sum_{x=0}^{k-1} q_k(x)f(x)\xi^x z^{\underline{\ell x}} \Big],
\end{array}\]
which leads to
\[\mathcal{R}(\xi,\ell;z)F(\xi,\ell;z) = \displaystyle\sum\limits_{k=0}^{r} \sum\limits_{x=0}^{k-1} p_k(x-k)f(x) \xi^{x+r-k} z^{\underline{\ell(x+r-k)}},\]
where $\mathcal R(\xi; \ell; z; \theta;\rho) = \displaystyle\sum_{k=0}^{r} \xi^{-k}z^{\underline{\ell(r-k)}}\rho^{\ell(r-k)} q_k(\theta)$. 
It follows directly that the operator $\mathcal{R}$ admits a representation matching the left-hand side of \eqref{discrete}. Applying Lemma~\ref{lemma3} causes the right-hand side of the equation to vanish.

Conversely, assume that the sequence $f$ is $p$-recursive. Transform the left-hand side of \eqref{discrete} as
\[\begin{array}{ll}
        \displaystyle\sum\limits_{j=0}^{r} P_j(\xi,\ell;z) \theta^j F(\xi,\ell;z) &=\displaystyle\sum\limits_{j=0}^{r} \sum\limits_{i=0}^{k_j} d_{i,j} \xi^{i} z^{\underline{i \ell}} \rho^{i \ell} \theta^j \sum\limits_{x=0}^\infty f(x) \xi^x z^{\underline{\ell x}} \\
        &= \displaystyle\sum\limits_{j=0}^{r} \sum\limits_{i=0}^{k_j} \sum\limits_{x=0}^\infty d_{i,j} x^j f(x) \xi^{x + i} z^{\underline{\ell( x + i)}} \\
        &=  \displaystyle\sum\limits_{j=0}^{r} \sum\limits_{i=0}^{k_j}  \sum\limits_{x=i}^\infty d_{i,j} (x-i)^j f(x-i) \xi^{x} z^{\underline{\ell x}}.
\end{array}\]
Note that the double sums $\displaystyle\sum_{i=0}^{k_j}\displaystyle\sum_{x=i}^{\infty}=\displaystyle\sum_{x=1}^{\infty} \displaystyle\sum_{i=0}^{\text{min}\{k_j,x-1\}}$, thus we obtain
\[ \displaystyle\sum_{j=0}^r P_j(\xi,\ell;z)\theta^j F(\xi,\ell;z)=\displaystyle\sum\limits_{x=1}^{\infty} \xi^x z^{\underline{\ell x}} \sum_{j=0}^r \sum_{i=0}^{\min\{k_j,x-1\}} d_{ij}(x-i)^j f(x-i)\]
Equating the coefficients of $z^{\underline{\ell x}}$ to zero results in  a difference equation for $f(x)$ of the form \eqref{difeqpol}, completing the proof.
\end{proof}

Consequently, we see that \eqref{discrete} serves as a discrete counterpart of \eqref{analog} for discrete generating series $F(\xi, \ell; z)$, which are accordingly called $D$-finite discrete generating series. It is worth noting that the study of $D$-finite generating series has been extended to multiple variables in \cite{Lipshits1989}, with conditions for their $D$-finiteness provided in \cite{Nekrasova2014, LeinartasNekrasova2016}. Additionally, the $D$-finiteness of sections of such series has been explored in \cite{LyapinCuchta2022, LyapinAkhtamova2021, Luzon2010}.

\section{Discrete Generating Series: Multidimensional Case}\label{sec:multidim}

An approach to constructing the general theory of discrete generating series and its connection with linear difference equations was introduced in \cite{Alekseev2023a, Akhtamova2024}.  Specifically, we define a discrete generating series for functions \( f\colon \mathbb{Z}^n \to \mathbb{C} \) and derive functional relations for such series.

Let $\mathbb{Z}_\geqslant$ denote the nonnegative integers, $\mathbb{Z}^n = \mathbb{Z} \times \cdots \times \mathbb{Z}$ be the $n$-dimensional integers, $\mathbb{Z}^n_\geqslant = \mathbb{Z}_\geqslant \times \cdots \times \mathbb{Z}_\geqslant$ for $n \in \mathbb{Z}_\geqslant$ be its nonnegative orthant. Throughout, we will use the multidimensional notation for convience of expressions: $x = (x_1, \ldots, x_n)\in \mathbb{Z}^n_\geqslant$, $z = (z_1, \ldots, z_n) \in \mathbb C^n$, $\xi = (\xi_1, \ldots, \xi_n) \in \mathbb C^n$, $\xi^x = \xi_1^{x_1} \cdots \xi_n^{x_n}$, $z^{\underline{x}} = z_1^{\underline{x_1}} \cdots z_n^{\underline{x_n}}$, $\ell = (\ell_1, \ldots, \ell_n) \in \mathbb Z^n_\geqslant$, $x! = x_1! \ldots x_n!$. We also will use $x \leqslant y$ for $x,y \in \mathbb{Z}^n$ componentwise, i.e. that $x_i \leqslant y_i$ for all $i=1,\ldots,n$.

Given a function $f \colon \mathbb{Z}^n_\geqslant \rightarrow \mathbb{C}$, we define the associated multidimensional discrete generating series of $f$ as
\[F(\xi; \ell; z) = \sum_{x_1=0}^{\infty} \ldots \sum_{x_n=0}^{\infty} f(x_1, \ldots, x_n) \xi_1^{x_1} \cdots \xi_n^{x_n} z_1^{\underline{\ell_1 x_1}} \cdots z_n^{\underline{\ell_n x_n}}.\]

{Let $p_{\alpha} \in \mathbb{C}[z]$ denote polynomials with complex coefficients.} The difference equation under consideration in this work is
\begin{equation}\label{equation}
    \sum_{\alpha \in A} p_\alpha(x) f(x + \alpha) = 0,
\end{equation}
where {set $A \subset\mathbb{Z}^n_\geqslant$ is finite and there is $m \in A$ such that for all $\alpha \in A$ the inequality $\alpha \leqslant m,$ which means $\alpha_j \leqslant m_j, j=1, \ldots, n,$ holds. Sometimes we will use an equivalent notation $0 \leqslant \alpha \leqslant m$ assuming that for some $\alpha$ coefficients $p_\alpha(x)$ vanish and only $p_m(x) \not\equiv 0$. 

We equip \eqref{equation} with initial data on a set named $X_m$ which is used often enough we introduce the notation $\mathbb{Z}_{\not\geqslant}$ as $X_m=\mathbb Z^n_\geqslant \setminus \left( m + \mathbb Z^n_\geqslant \right) = \left\{ x \in \mathbb{Z}^n_{\geqslant} \colon x \not\geqslant m\right\}$ and we define the initial data function $\varphi: X_m \to \mathbb C$ so that
\begin{equation}
f(x)=\varphi(x), \quad x \in X_m.\label{initialdata}
\end{equation}

\begin{figure}
\centering

\begin{subfigure}{0.3\textwidth}
\begin{tikzpicture}
\draw[pattern=north west lines, pattern color=black] (1,1) rectangle (3,2);
\draw (0,0) -- (0,2);
\draw (0,0) -- (3,0);
\node at (1,1) {\textbullet};
\node at (1.2,0.85) {$m$};
\draw[dashed] (0,1)--(3,1);
\draw[dashed] (1,0)--(1,2);
\draw[color=white,thick] (1,2)--(3,2);
\draw[color=white,thick] (3,1)--(3,2);
\end{tikzpicture}
\caption{$x \geqslant m$}
\end{subfigure}
\begin{subfigure}{0.3\textwidth}
\begin{tikzpicture}
\draw (0,0) -- (0,2);
\draw (0,0) -- (3,0);
\node at (1,1) {\textbullet};
\node at (1.2,0.85) {$m$};
\draw[pattern=north west lines, pattern color=black] (0,0) rectangle (1,1);
\draw[dashed] (0,1)--(3,1);
\draw[dashed] (1,0)--(1,2);
\end{tikzpicture}
\caption{$x \leqslant m$}
\end{subfigure}
\begin{subfigure}{0.3\textwidth}
\begin{tikzpicture}
\node at (1,1) {\textbullet};
\node at (1.2,1.2) {$m$};
\draw[pattern=north east lines, pattern color=black] (0,0)--(3,0)--(3,1)--(1,1)--(1,2)--(0,2)--(0,0);
\draw[color=white] (1,1)--(1,2);
\draw[color=white] (1,1)--(3,1);
\node at (1,1) {\textbullet};
\draw[dashed] (1,1)--(1,2);
\draw[dashed] (1,1)--(3,1);
\draw[color=white,thick] (0,2)--(1,2);
\draw[color=white,thick] (3,0)--(3,1);
\end{tikzpicture}
\caption{$x \not\geqslant m$}
\end{subfigure}
\caption{Illustration of the sets $x \geqslant m$, $x \leqslant m$, and $x\not\geqslant m$.}
\end{figure}

By analogy with the one-dimensional case, the Cauchy problem \cite{Nekrasova2015} involves finding a solution to the difference equation \eqref{equation} that coincides with $\varphi$ on $X_m$, i.e., $f(x) = \varphi(x)$ for all $x \in X_m$.

\subsection{Discrete generating series for linear difference equations with constant coefficients in higher dimension}\label{Section2}

In this section, we consider higher-order homogeneous difference equations with constant coefficients
\begin{equation}\label{constantcoeff}
    \sum_{\alpha \in A} c_{\alpha} f(x + \alpha) = 0,
\end{equation}
which generalize \eqref{difeqconst}.
%
The shift operator is
\begin{equation*}\label{shiftopscriptp}
\mathcal P(\xi,\ell;z) = \sum_{0 \leqslant \alpha \leqslant m} c_\alpha \xi^\alpha z^{\underline{\ell \alpha}} \rho^{\ell \alpha},
\end{equation*}
the truncation for $\tau \in \mathbb Z^n$ is defined by the formula
\[\mathcal P_\tau(\xi,\ell;z) = \sum_{\substack{0 \leqslant \alpha \leqslant m \\ \alpha \ngeqslant
    \tau}} c_\alpha \xi^\alpha z^{\underline{\ell \alpha}} \rho^{\ell \alpha},\]
and the discrete generating series of the initial data for $\tau \in X_m$ is
\begin{equation*}
\Phi_\tau(\xi, \ell; z) = \sum_{x \ngeqslant \tau} \varphi(x) \xi^x z^{\underline{\ell x}}.\label{genfuncinitdata}
\end{equation*}
Let $\delta_j\colon x \to x + e^j$ be the forward shift operator for $j = 1, \ldots, n$ with multidimensional notation $\delta^\alpha = \delta_1^{\alpha_1} \ldots \delta_n^{\alpha_n}$ and define the polynomial difference operator
\[P(\delta) = \sum_{0 \leqslant \alpha \leqslant m} c_\alpha \delta^{\alpha}.\]
With this notation, equation \eqref{constantcoeff} is represented compactly as
\[P(\delta^{-1}) f(x) = 0, \ \ \ x \geqslant m.\]
The study of generating series of the form \(\sum_x f(x) z^x\) has been extensively explored for both single and multiple variables. One of the earliest and most practical formulas for deriving such generating series, utilizing the characteristic polynomial and the initial data function, was established in \cite{Leinartas2009}. Additionally, several effective computer algebra techniques for this purpose have been developed and implemented in \cite{KytmanovLyapinSadykov2017, ApanovichLyapinShadrin2021}. In this work, we will establish analogous formulas for the discrete generating series \(F(\xi,\ell; z)\).

\begin{theorem}
    The discrete generating series $F(\xi, \ell; z)$ for the solution to the Cauchy problem for equation \eqref{constantcoeff} with initial data \eqref{initialdata} satisfies the functional equations
    \begin{align*}
          \mathcal P(\xi,\ell;z) F(\xi, \ell; z)
          &= \sum_{0 \leqslant \alpha \leqslant m} c_\alpha \xi^\alpha z^{\underline{\ell \alpha}} \rho^{\ell \alpha} \Phi_{m-\alpha}(\xi, \ell; z)  \\
          &= \sum_{x \ngeqslant m} P(\delta^{-1}) \varphi(x)  \xi^x z^{\underline{\ell x}}  \\
          &= \sum_{x \ngeqslant m} \mathcal P_{m - x}(\xi,\ell;z) \varphi(x) z^{\underline{\ell x}}. 
    \end{align*}
\end{theorem}

\begin{proof}
  Multiplying \eqref{constantcoeff} by $\xi^x z^{\underline{\ell x}}$ and summing over $x \geqslant m$, we obtain
  \begin{align*}
      0&=\sum_{x \geqslant m} \sum_{0 \leqslant \alpha \leqslant m} c_\alpha f(x-\alpha) \xi^x z^{\underline{\ell x}}  \\
    &=\sum_{0 \leqslant \alpha \leqslant m} c_\alpha \sum_{x \geqslant m}  f(x-\alpha) \xi^x z^{\underline{\ell x}}.
\end{align*}
Substituting $x$ with $x+\alpha$ and some routine manipulations complete the proof \cite{Akhtamova2024}.
\end{proof}

For $z = (z_1, \ldots, z_n)$ we denote the projection operator $$\pi_j z = (z_1, \ldots, z_{j-1}, 0, z_{j+1}, \ldots, z_n)$$ 
the projection of the discrete generating series
\[\pi_j F(\xi; \ell; z) = F(\xi; \ell; \pi_j z) = \sum_{\substack{x \geqslant 0\\ x_j = 0}} f(x) \xi^x z^{\underline{\ell x}},\]
and the combined projection $\Pi = (1-\pi_1) \circ \cdots \circ (1-\pi_n)$ as the composition of $1-\pi_j$ for all $j = 1, \ldots, n$.

For the next result, we introduce the symbols $I=(1,1,\ldots,1)$ and the unit vectors $e_j=(0,\ldots,0,1,0,\ldots, 0)$ for $j=1,2,\ldots,n$ which is nonzero only the $j$th component. Useful properties of the combined projection $\Pi$ were derived in \cite{Akhtamova2024}:
\begin{align*}
   \Pi \sum_{x \geqslant 0} f(x) \xi^x z^{\underline{\ell x}} = \sum_{x \geqslant I} f(x) \xi^x z^{\underline{\ell x}},
\end{align*}
and
\begin{align*}
        \Pi \xi_j z_j^{\underline{\ell_j}} \rho^{\ell_j} F(\xi; \ell; z) = \sum_{x \geqslant I} f(x-e^j) \xi^x z^{\underline{\ell x}}.
\end{align*}

{We introduce the inner  product $$\langle c, \xi z^{\underline{\ell}} \rho^\ell \rangle = c_1 \xi_1 z_1^{\underline{\ell_1}} \rho_1^{\ell_1} + \cdots + c_n \xi_n z_n^{\underline{\ell_n}} \rho_n^{\ell_n}$$ and
\[\langle c, \delta^{-I} \rangle = c_1 \delta_1^{-1} + \cdots + c_n \delta_n^{-1}.\] }
We are now prepared to prove an analogue of \cite[Theorem~1.1]{LyapinChandragiri2019}.
\begin{theorem}
    The following formula holds:
    \begin{equation*}
        \Pi \left[ (1 - \langle c, \xi z^{\underline{\ell}} \rho^\ell \rangle) F(\xi;\ell;z) \right] =
        \sum_{x \geqslant I} (1- \langle c, \delta^{-I}\rangle ) f(x) \xi^x z^{\underline{\ell x}}.
    \end{equation*}
\end{theorem}

\begin{proof}
Applying $\Pi$ to $(1 - \langle c, \xi z^{\underline{\ell}} \rho^\ell \rangle) F(\xi;\ell;z)$ yields the proof \cite{Akhtamova2024}.
\end{proof}

The following corollary is straightforward.
\begin{corollary}
If $f$ solves $(1- \langle c, \delta^{-I}\rangle ) f(x)=0$, then
    \begin{equation*}
        \Pi \left[ (1 - \langle c, \zeta \rangle) F(\xi;\ell;z) \right] = 0.
    \end{equation*}
\end{corollary}

\subsection{Discrete generating series for linear difference equations with polynomial coefficients}\label{sec3}

We define the forward partial difference operators $\Delta_j$ by
\[\Delta_j F(z) = F(z+e^j) - F(z), \qquad j = 1, \ldots, n.\]
If $z^{\underline{x}} = z_1^{\underline{x_1}}\ldots z_n^{\underline{x_n}}$, then in terms of the multi-dimensional notation, we have $\Delta_j z^{\underline{x}} = x_j z^{\underline{x-e^j}}$. Now compute
\begin{align*}
    \Delta_j F(\xi; \ell; z) &=\Delta_j \sum_{x \geqslant 0} f(x) \xi^x z^{\underline{\ell x}} \\
    &= \sum_{x \geqslant 0} f(x) \xi^x \Delta_j z^{\underline{\ell x}} \\
    &= \sum_{x \geqslant 0} \ell_j x_j f(x) \xi^x  z^{\underline{\ell x - e^j}}.
\end{align*}

We denote by $\rho_j$ the componentwise backwards jump in the $j$th component, i.e. $\rho_j F(z) = F(z - e^j)$, and we define the componentwise operators $\theta_j = \ell_j^{-1} z_j \rho_j \Delta_j$, which generalizes the single-variable one defined earlier in \eqref{thetaofdiscretemonomial}. Useful properties for $p(\theta)=\displaystyle\sum_{\alpha \in A \subset \mathbb{Z}_{\geq}^n} c_{\alpha}\theta^{\alpha}$ were proven in \cite{Akhtamova2024}: if $k = (k_1, \ldots, k_n)\in \mathbb{Z}^n_{\geqslant}$, then 
    \begin{align*}
        \theta^k F(\xi; \ell; z) = \sum_{x \geqslant 0} x^{k} f(x) \xi^x z^{\underline{\ell x}},
    \end{align*}
and
    \begin{align*}
        p(\theta) F(\xi; \ell; z) = \sum_{x \geqslant 0} p(x) f(x) \xi^x z^{\underline{\ell x}}.
    \end{align*}
We define for $q_{\alpha}(x)=p_{\alpha}(x+\alpha)$ an operator $\mathcal{R}_A$ by
$$\mathcal R_A(\xi, \ell; z) = \sum_{\alpha \in A} q_{\alpha}(\theta) \xi^\alpha  z^{\underline{\ell \alpha}} \rho^{\ell \alpha}.$$
This framework leads to the following theorem.
\begin{theorem}
    The discrete generating series $F(\xi,\ell;z)$ of the Cauchy problem for equation \eqref{equation} with initial data \eqref{initialdata} satisfies the functional equation
    \begin{equation*}
          \mathcal R_A(\xi, \ell; z) F(\xi, \ell; z) = \sum_{\alpha \in A} \sum_{x \ngeqslant m-\alpha}  p_\alpha(x - \alpha) \varphi(x) \xi^{x} z^{\underline{\ell(x+\alpha)}}.
      \end{equation*}
\end{theorem}
\begin{proof}
      Similar to the proof of Theorem~\ref{thm1}, we multiply \eqref{equation} by $\xi^x z^{\underline{\ell x}}$ and sum over $x \geqslant m$,
      replacing $x$ with $x+\alpha$ then routine algebraic manipulation completes the proof.
\end{proof}

\section{Examples}\label{sec:examples}

\begin{example}
Consider the tribonacci equation
\begin{align*}
f(x+3) - f(x+2) - f(x+1) - f(x) = 0,
\end{align*}
multiply both sides by $\xi^{x+3} z^{\underline{\ell (x+3)}}$, sum over $x \geqslant 0$, and reindex to obtain
\begin{align*}
   0&=\sum_{x \geqslant 3}^\infty f(x)\xi^{x} z^{\underline{\ell x}}  - \xi z^{\underline{\ell}}\rho^{\ell}\sum_{x \geqslant 2}^\infty f(x)\xi^{x} z^{\underline{\ell x}} \\ &\phantom{=}-\xi^2 z^{\underline{2 \ell}}\rho^{2 \ell}\sum_{x \geqslant 1}^\infty f(x)\xi^{x} z^{\underline{\ell x}} -  \xi^3 z^{\underline{3 \ell}}\rho^{3 \ell}\sum_{x \geqslant 0}^\infty f(x)\xi^{x} z^{\underline{\ell x}}.
\end{align*}
Using the notation of $F(\xi,\ell;z)$ yields
\begin{align*}
    0&=F(\xi,\ell;z)-f(0)-f(1)\xi z^{\underline{\ell}}-f(2)\xi^2 z^{\underline{2 \ell}} \\ &\phantom{=}-\xi z^{\underline{\ell}} \rho^{\ell} \Big(F(\xi,\ell;z)-f(0)-f(1)\xi z^{\underline{\ell}}\Big) \\ &\phantom{=}-\xi^2 z^{\underline{2 \ell}}\rho^{2 \ell} \Big(F(\xi,\ell;z)-f(0)\Big) \\ &\phantom{=}-\xi^3 z^{\underline{3 \ell}}\rho^{3 \ell}F(\xi,\ell;z).
\end{align*}
Finally, the functional equation for the discrete generating series for tribonacci numbers is
\begin{align*}
     \bigg( 1-\xi z^{\underline{\ell}}\rho^{\ell} -\xi^2 z^{\underline{2 \ell}}\rho^{2 \ell} -\xi^3 z^{\underline{3 \ell}}\rho^{3 \ell} \bigg)  F(\xi;\ell;z) &=f(0)+f(1)\xi z^{\underline{\ell}}+f(2)\xi^2 z^{\underline{2 \ell}} \\  &\phantom{=}-\xi z^{\underline{\ell}} \rho^{\ell}(f(0)+f(1)\xi z^{\underline{\ell}})  \\ &\phantom{=}-\xi^2 z^{\underline{2 \ell}}\rho^{2 \ell}f(0) .
\end{align*}
\end{example}

This example can be naturally extended to broader classes of generalized Fibonacci numbers, including those defined by higher-order recurrences, weighted sums, or non-linear terms. Such generalizations have been extensively studied, for instance, in the works of \cite{Kizilates2024} and \cite{Kizilates2025}.

\begin{example}[The second problem of Schr\"{o}der \cite{schroeder1870}]
Let $f(x)$ be the total number of ways of dividing an $(x+1)$-gon  by diagonals not intersecting in the interior. Function $f(x)$ solves the difference equation \cite[p. 57]{Comtet1974}
\begin{align*}
(x+2)f(x+2) - 3(2x+1)f(x+1) + (x-1)f(x) = 0.
\end{align*}
Multiplying  both sides by $\xi^{x+2} z^{\underline{\ell (x+2)}}$ and summing over $x \geqslant 0$ yield
\begin{align*}
    0 =&\sum\limits_{x=0}^\infty (x+2)f(x+2) \xi^{x+2} z^{\underline{\ell (x+2)}} \\ - &\sum\limits_{x=0}^\infty 3(2x+1) f(x+1) \xi^{x+2} z^{\underline{\ell (x+2)}} \\ + &\sum\limits_{x=0}^\infty (x-1) f(x) \xi^{x+2} z^{\underline{\ell (x+2)}}.
\end{align*}
Reindexing and some routine manipulations yields
\begin{multline*}
    \displaystyle\sum\limits_{x=2}^\infty x f(x) \xi^{x} z^{\underline{\ell x}}- \xi z^{\underline{\ell}}\rho^\ell \displaystyle\sum\limits_{x=1}^\infty 3(2x-1) f(x) \xi^{x} z^{\underline{\ell x}} + \xi^2 z^{\underline{2 \ell}} \rho^{2 \ell} \displaystyle\sum\limits_{x=0}^\infty (x-1) f(x) \xi^{x} z^{\underline{\ell x}} \\  
    =\theta\Big(F(\xi;\ell;z)-f(0)-f(1)\xi z^{\underline{\ell}}\Big) -3(2\theta-1)\xi z^{\underline{\ell}}\rho^\ell \Big(F(\xi;\ell;z)-f(0)\Big)\\ 
    + (\theta-1)\xi^2 z^{\underline{2\ell}}\rho^{2\ell} F(\xi;\ell;z).
\end{multline*}
This results in the functional relation
\begin{multline*}
(\theta - 3(2 \theta - 1)\xi z^{\underline{\ell}} \rho^{\ell} + (\theta - 1)\xi^2 z^{\underline{2 \ell}}\rho^{2 \ell})F(\xi;\ell;z) \\ =\theta(f(0) + f(1)\xi z^{\underline{\ell}}) - 3(2\theta-1)f(0).
\end{multline*}
By Lemma~\ref{lemma3}, applying $\theta(\theta-1)$ to each side yields
\begin{align*}
(\theta-1)\theta \!\Big[ \theta - 3(2 \theta - 1)\xi z^{\underline{\ell}} \rho^{\ell} + (\theta - 1)\xi^2 z^{\underline{2 \ell}}\rho^{2 \ell} \Big] F(\xi;\ell;z)  
= 0.
\end{align*}
\end{example}

\begin{example}
We consider a functional equation
\begin{equation*}
\begin{gathered}
\xi^2 z^{\underline{2\ell}}\rho^{2\ell}\theta^3F(\xi;\ell;z) + \xi^3 z^{\underline{3\ell}}\rho^{3\ell}\theta F(\xi;\ell;z) = 0.
\end{gathered}
\end{equation*}
Substitution $F(\xi;\ell;z)$ and some simple manipulations yields
\begin{equation*}
\sum_{x \geqslant 0}x^3f(x)\xi^{x+2}z^{\underline{\ell (x+2)}} + \sum_{x \geqslant 0}xf(x)\xi^{x+3}z^{\underline{\ell (x+3)}} = 0.
\end{equation*}
Changing the summation indices in each sum yields
\begin{equation*}
\sum_{x \geqslant 2}(x-2)^3f(x-2)\xi^{x}z^{\underline{\ell x}} + \sum_{x \geqslant 3}(x-3)f(x-3)\xi^{x}z^{\underline{\ell x}} = 0.
\end{equation*}
Since the first sum equals 0 for $x=2$, we can write the entire expression as a single sum for $x \geqslant 3$ as
\begin{equation*}
\sum_{x \geqslant 3} \left( (x-2)^3f(x-2)+ (x-3)f(x-3)\right) \xi^{x}z^{\underline{\ell x}} = 0,
\end{equation*}
that leads to a difference equation
\begin{equation*}
(x-2)^3f(x-2) + (x-3)f(x-3) = 0,
\end{equation*}
or equivalently
\begin{equation*}
\begin{gathered}
(x+1)^3f(x+1) + xf(x) = 0.
\end{gathered}
\end{equation*}
\end{example}

\begin{example}
\label{example4} We will derive the functional equation for the discrete generating series $F(;;z_1,z_2):=F(\xi_1, \xi_2; \ell_1, \ell_2; z_1, z_2)$ for the basic combinatorial recurrence
\begin{equation}\label{basicrecurrence}
    f(x_1+1, x_2+1) - f(x_1, x_2+1) - f(x_1+1, x_2) = 0.
\end{equation}
Multiplying both sides of \eqref{basicrecurrence} by $\xi_1^{x_1+1} \xi_2^{x_2+1} z_1^{\underline{\ell_1 (x_1+1)}} z_2^{\underline{\ell_2 (x_2+1)}}$ and summing over $(x_1, x_2) \geqslant (0,0)$ yields
\begin{align*}
    0=&\sum_{(x_1, x_2) \geqslant (0,0)} f(x_1+1, x_2+1)\xi_1^{x_1+1} \xi_2^{x_2+1} z_1^{\underline{\ell_1 (x_1+1)}} z_2^{\underline{\ell_2 (x_2+1)}} \\ - &\sum_{(x_1, x_2) \geqslant (0,0)}f(x_1 , x_2+1)\xi_1^{x_1+1} \xi_2^{x_2+1} z_1^{\underline{\ell_1 (x_1+1)}} z_2^{\underline{\ell_2 (x_2+1)}} \\- &\sum_{(x_1, x_2) \geqslant (0,0)}f(x_1+1, x_2)\xi_1^{x_1+1} \xi_2^{x_2+1} z_1^{\underline{\ell_1 (x_1+1)}} z_2^{\underline{\ell_2 (x_2+1)}}.
\end{align*}
We consider each sum separately: first,
\[\begin{array}{ll}
    \displaystyle\sum_{(x_1, x_2) \geqslant (0,0)} f(x_1+1, x_2+1)\xi_1^{x_1+1} \xi_2^{x_2+1} z_1^{\underline{\ell_1 (x_1+1)}} z_2^{\underline{\ell_2 (x_2+1)}} \\
     \phantom{\quad}=\displaystyle\sum_{(x_1, x_2) \geqslant (1,1)} f(x_1, x_2)\xi_1^{x_1} \xi_2^{x_2} z_1^{\underline{\ell_1 x_1}} z_2^{\underline{\ell_2 x_2}} \\
    \phantom{\quad}= F(;; z_1, z_2) - F(;; 0, z_2) - F(;; z_1, 0) + F(;; 0, 0),\\
\end{array}\]
second,
\[\begin{array}{ll}
\displaystyle\sum_{(x_1, x_2) \geqslant (0,0)}f(x_1 , x_2+1)\xi_1^{x_1+1} \xi_2^{x_2+1} z_1^{\underline{\ell_1 (x_1+1)}} z_2^{\underline{\ell_2 (x_2+1)}} \\
    \phantom{\quad}=\displaystyle\sum_{(x_1, x_2) \geqslant (0,1)}f(x_1, x_2)\xi_1^{x_1+1} \xi_2^{x_2} z_1^{\underline{\ell_1 (x_1+1)}} z_2^{\underline{\ell_2 x_2}} \\
    \phantom{\quad}=\xi_1 z_1^{\underline{\ell_1}} \rho_1^{\ell_1} \displaystyle\sum_{(x_1, x_2) \geqslant (0,1)}f(x_1, x_2)\xi_1^{x_1} \xi_2^{x_2} z_1^{\underline{\ell_1 x_1}} z_2^{\underline{\ell_2 x_2}} \\
    \phantom{\quad}=\xi_1 z_1^{\underline{\ell_1}} \rho_1^{\ell_1} \big( F(;; z_1, z_2) - F(;; z_1, 0)\big),\\
\end{array}\]
and third,
\[\begin{array}{ll}
\displaystyle\sum_{(x_1, x_2) \geqslant (0,0)}f(x_1+1, x_2)\xi_1^{x_1+1} \xi_2^{x_2+1} z_1^{\underline{\ell_1 (x_1+1)}} z_2^{\underline{\ell_2 (x_2+1)}} \\
    \phantom{\quad}=\displaystyle\sum_{(x_1, x_2) \geqslant (1,0)}f(x_1, x_2)\xi_1^{x_1} \xi_2^{x_2+1} z_1^{\underline{\ell_1 x_1}} z_2^{\underline{\ell_2 (x_2+1)}} \\
    \phantom{\quad}=\xi_2 z_2^{\underline{\ell_2}} \rho_2^{\ell_2} \displaystyle\sum_{(x_1, x_2) \geqslant (1,0)}f(x_1, x_2)\xi_1^{x_1} \xi_2^{x_2} z_1^{\underline{\ell_1 x_1}} z_2^{\underline{\ell_2 x_2}} \\
    \quad =\xi_2 z_2^{\underline{\ell_2}} \rho_2^{\ell_2} \big( F(;; z_1, z_2) - F(;; 0, z_2)\big).
\end{array}\]
Finally, we get
\begin{equation*}
\begin{gathered}
    F(;; z_1, z_2) - F(;; 0, z_2) - F(;; z_1, 0) + F(;; 0, 0) \\ - \xi_1 z_1^{\underline{\ell_1}} \rho_1^{\ell_1} \big( F(;; z_1, z_2)  - F(;; z_1, 0)\big)  - \xi_2 z_2^{\underline{\ell_2}} \rho_2^{\ell_2} \big( F(;; z_1, z_2) - F(;; 0, z_2)\big) = 0,
\end{gathered}
\end{equation*}
which yields the following functional equation on $F(;;z_1, z_2)$:
\begin{equation*}
\begin{gathered}
    (1 - \xi_1 z_1^{\underline{\ell_1}} \rho_1^{\ell_1} - \xi_2 z_2^{\underline{\ell_2}} \rho_2^{\ell_2}) F(;; z_1, z_2) \\
    -(1 - \xi_2 z_2^{\underline{\ell_2}} \rho_2^{\ell_2}) F(;; 0, z_2) - (1 - \xi_1 z_1^{\underline{\ell_1}} \rho_1^{\ell_1}) F(;; z_1, 0) + F(;; 0, 0) = 0.
\end{gathered}
\end{equation*}
\end{example}

\section{Conclusion}\label{sec:conclusion}
We have laid the foundation for the theory of discrete generating series in the context of multidimensional difference equations with polynomial coefficients. By introducing a multidimensional polynomial shift operator, we derived three functional equations that these discrete generating series satisfy, uncovering key structural properties. The falling factorial functions in this framework are referred to as generalized polynomials in the context of time scales. This connection hints at potential extensions of this research to time scales, a direction that was arguably foreshadowed by the definition of moment generating series for distributions in \cite{matthewsthesis}. A particularly intriguing question is identifying the appropriate analogue of \eqref{equation} for an arbitrary time scale. Exploring this problem from a generating series perspective may yield new insights and deepen our understanding of the subject.

%
%

\end{document}